\documentclass[reqno,10pt]{amsart}

\newtheorem{thm}{Theorem}[section]

\newtheorem{lem}[thm]{Lemma}

\newtheorem{rek}[thm]{Remark}

\newcommand\be{\begin{equation}}
\newcommand\ee{\end{equation}}
\newcommand\bea{\begin{eqnarray}}
\newcommand\eea{\end{eqnarray}}
\newcommand\bi{\begin{itemize}}
\newcommand\ei{\end{itemize}}
\newcommand\ben{\begin{enumerate}}
\newcommand\een{\end{enumerate}}
\newcommand\bc{\begin{center}}
\newcommand\ec{\end{center}}
\newcommand\ba{\begin{array}}
\newcommand\ea{\end{array}}

\newcommand{\Qoft}{\Bbb{Q}(T)}  %use in linux
\newcommand{\Z}{\ensuremath{\Bbb{Z}}}
\newcommand{\Q}{\Bbb{Q}}

\newcommand{\ga}{\alpha}     %lowercase alpha
\newcommand{\lap}{\triangle}    %laplacian
\newcommand{\hphi}{\widehat{\phi}}  %phi^
\newcommand{\phiint}{\int_{-\infty}^{\infty} \phi(y) dy}

\newcommand{\twocase}[5]{#1 \begin{cases} #2 & \text{#3}\\ #4
&\text{#5} \end{cases}   }

\newcommand{\foh}{\frac{1}{2}}  %onehalf
\newcommand{\js}[1]{ { \underline{#1} \choose p} }
\newcommand{\zsum}[1]{ \sum_{#1 = 0}^{p-1} }
\newcommand{\osum}[1]{ \sum_{#1 = 1}^{p-1} }

\renewcommand{\theequation}{\thesection.\arabic{equation}}

\begin{document}

\title{Constructing
Elliptic Curves over $\Q(T)$ with Moderate Rank}

\author{Scott Arms}
\email{arms@math.ohio-state.edu}
\address{Department of Mathematics, The Ohio State University, Columbus, OH, $43210$, U.S.A.}

\author{Steven J. Miller}
 \email{sjmiller@math.ohio-state.edu}
 \address{Department of Mathematics, The Ohio State University, Columbus, Ohio, 43210, U.S.A.}

\author{\'{A}lvaro Lozano-Robledo}
\email{alozano@math.bu.edu}
\address{Department of Mathematics, Boston University, Boston, MA, U.S.A.}

%\classification{11M26, 11G05, 11G40, 11M26.}

\date{}

\begin{abstract}
We give several new constructions for moderate rank elliptic
curves over $\Q(T)$. In particular we construct infinitely many
rational elliptic surfaces (not in Weierstrass form) of rank $6$
over $\Q$ using polynomials of degree two in $T$. While our method
generates linearly independent points, we are able to show the
rank is exactly $6$ \emph{without} having to verify the points are
independent. The method generalizes; however, the higher rank
surfaces are not rational, and we need to check that the
constructed points are linearly independent.
\end{abstract}

\maketitle

%\tableofcontents

%!!!!!!!!!!!!!!!!!!!!!!!!!!!!!!!!!!!!!!!new section!!!!!!!!!!!!!!!!!!!!!!!!!!!!!!!!
%!!!!!!!!!!!!!!!!!!!!!!!!!!!!!!!!!!!!!!!new section!!!!!!!!!!!!!!!!!!!!!!!!!!!!!!!!
%!!!!!!!!!!!!!!!!!!!!!!!!!!!!!!!!!!!!!!!new section!!!!!!!!!!!!!!!!!!!!!!!!!!!!!!!!

\section{Introduction}
\setcounter{equation}{0}

Consider the elliptic curve $\mathcal{E}$ over $\Q(T)$:
\be\label{eq:EoverQt} y^2 + a_1(T)xy + a_3(T)y\ =\ x^3 + a_2(T)x^2
+ a_4(T)x + a_6(T), \ee where $a_i(T) \in \Z[T]$. By evaluating
these polynomials at integers, we obtain elliptic curves over
$\Q$. By Silverman's Specialization Theorem, for large $t\in\Z$
the Mordell-Weil rank of the fiber $\mathcal{E}_t$ over $\Q$ is at
least that of the curve $\mathcal{E}$ over $\Qoft$.

%\subsection{Other Methods to Construct Families with Rank}
For comparison purposes, we briefly describe other methods to
construct curves with rank. Mestre \cite{Mes1,Mes2} considers a
$6$-tuple of integers $a_i$ and defines $q(x) = \prod_{i=1}^6 (x -
a_i)$ and $p(x,T) = q(x-T)q(x+T)$. There exist polynomials
$g(x,T)$ of degree $6$ in $x$ and $r(x,T)$ of degree at most $5$
in $x$ such that $p(x,T) = g^2(x,T) - r(x,T)$. Consider the curve
$y^2 = r(x,T)$ over $\Q(T)$. If $r(x,T)$ is of degree $3$ or $4$
in $x$, we obtain an elliptic curve with points $P_{\pm i}(T) =
\left(\pm T + a_i, g(\pm T + a_i)\right)$. If $r(x,T)$ has degree
$4$  we may need to change variables to make the coefficient of
$x^4$ a perfect square (see \cite{Mor}, page 77). Two $6$-tuples
that work are $(-17,-16,10,11,14,17)$ and $(399,380,352,47,4,0)$
(see \cite{Na1}). Curves of rank up to $14$ over $\Q(T)$ have been
constructed this way, and using these methods Nagao \cite{Na1} has
found an elliptic curve of rank at least $21$ and Fermigier
\cite{Fe2} one of rank at least 22 over $\Q$. Shioda \cite{Sh}
gives explicit constructions for not only rational curves of rank
$2,4,6,7$ and $8$ over $\Q(T)$, but generators of the Mordell-Weil
groups.

We now describe the idea of our method. For $\mathcal{E}$ as in
\eqref{eq:EoverQt}, define
\begin{eqnarray}
A_\mathcal{E}(p) \ = \  \frac{1}{p} \sum_{t=0}^{p-1} a_t(p),
\end{eqnarray}
with $a_t(p) = p+1 - N_t(p)$, where $N_t(p)$ is the number of
points in $\mathcal{E}_t(\mathbb{F}_p)$ (we set $a_t(p)=0$ when
$p\mid \Delta(t)$). Rosen and Silverman \cite{RS} prove a version
of a conjecture of Nagao \cite{Na1} which relates
$A_{\mathcal{E}}(p)$ to the rank of $\mathcal{E}$ over $\Q(T)$.

\begin{thm}[Rosen-Silverman]\label{thm:RoSi}
Let $\mathcal{E}:y^2 = x^3 + A(T)x + B(T)$, $A,B \in \Z[T]$, and
assume Tate's conjecture (known if $\mathcal{E}$ is a rational
elliptic surface over $\Q$) for $\mathcal{E}$. Then
\begin{eqnarray}
\lim_{X \to \infty} \frac{1}{X} \sum_{p \leq X} -A_\mathcal{E}(p)
\log p \ = \  \text{\emph{rank}} \ \mathcal{E}(\Q(T)).
\end{eqnarray}
\end{thm}

An elliptic curve $\mathcal{E}$ (as in the previous theorem) is a
rational elliptic surface over $\Q$ if and only if one of the
following holds:
\begin{enumerate}
  \item $0 < \max\{3 \text{deg}\ A(T), 2\text{deg}\ B(T)\} < 12$.
  \item $3\text{deg}\ A(T) = 2\text{deg}\ B(T) = 12$ and
$\text{ord}_{T=0}T^{12} \Delta(T^{-1}) = 0$
\end{enumerate}
(see \cite{Mir,RS}). In this paper we construct special rational
elliptic surfaces where we are able to evaluate
$A_{\mathcal{E}}(p)$ exactly. For these surfaces, we have
$A_\mathcal{E}(p)$ $=$ $-r + O(\frac{1}{p})$. By Theorem
\ref{thm:RoSi} and the Prime Number Theorem, we can conclude that
the constant $r$ is the rank of $\mathcal{E}$ over $\Q(T)$.

The novelty of this approach is that by forcing $A_\mathcal{E}(p)$
to be essentially constant, provided $\mathcal{E}$ is a rational
elliptic surface over $\Q$, we can immediately calculate the
Mordell-Weil rank \emph{without} having to specialize points and
calculate height matrices. Further, we obtain an exact answer for
the rank, and not a lower bound.

If the degrees of the defining polynomials of $\mathcal{E}$ are
too large, our results are conditional on Tate's conjecture if we
are able to evaluate $A_{\mathcal{E}}(p)$. In many cases, however,
we are unable to evaluate $A_{\mathcal{E}}(p)$ to the needed
accuracy. Our method does generate candidate points, which upon
specialization yield lower bounds for the rank. In this manner,
curves of rank up to 8 over $\Q(T)$ have been found.

Modifications of our method may yield curves with higher rank over
$\Q(T)$, though to \emph{find} such curves requires solving very
intractable non-linear Diophantine equations and then specializing
the points and calculating the height matrices to see that they
are independent over $\Q(T)$.

For additional constructions, especially for lower rank curves
over $\Q(T)$, see \cite{Fe2}. For a good survey on ranks of
elliptic curves, see \cite{RuS}.

%!!!!!!!!!!!!!!!!!!!!!!!!!!!!!!!!!!!!!!!new section!!!!!!!!!!!!!!!!!!!!!!!!!!!!!!!!
%!!!!!!!!!!!!!!!!!!!!!!!!!!!!!!!!!!!!!!!new section!!!!!!!!!!!!!!!!!!!!!!!!!!!!!!!!
%!!!!!!!!!!!!!!!!!!!!!!!!!!!!!!!!!!!!!!!new section!!!!!!!!!!!!!!!!!!!!!!!!!!!!!!!!

%!!!!!!!!!!!!!!!!!!!!!!!!!!!!!!!!!!!!!!!new section!!!!!!!!!!!!!!!
%!!!!!!!!!!!!!!!!!!!!!!!!!!!!!!!!!!!!!!!new section!!!!!!!!!!!!!!!
%!!!!!!!!!!!!!!!!!!!!!!!!!!!!!!!!!!!!!!!new section!!!!!!!!!!!!!!!

\section{Constructing Rank $6$ Rational Surfaces over $\Q(T)$}
\setcounter{equation}{0}

\subsection{Idea of the Construction}

The main idea is as follows: we can explicitly evaluate linear and
quadratic Legendre sums; for cubic and higher sums, we  cannot in
general  explicitly evaluate the sums. Instead, we have bounds
(Hasse, Weil) exhibiting large cancellation.

The goal is to cook up curves $\mathcal{E}$ over $\Q(T)$ where we
have linear and quadratic expressions in $T$. We can evaluate
these expressions exactly by a standard lemma on Quadratic
Legendre Sums (see Lemma \ref{labquadlegsum} of the appendix for a
proof),
%STEVE:(\textbf{this is a straightforward, standard result: do we want to
%just quote it, or keep the paper mostly self-contained by
%including the proof? if we are aiming for a short paper, I'd say
%eliminate the proof}
 which states that if $a$ and $b$ are not
both zero mod $p$ and $p
> 2$, then for $t\in\Z$
\begin{equation}
\twocase{\zsum{t} \js{at^2 + bt + c} \ = \ }{(p-1)\js{a}}{if $p |
 (b^2 - 4ac)$}{-\js{a}}{otherwise.}
\end{equation} Thus if $p|(b^2-4ac)$, the summands are $\js{a(t-t')^2} = \js{a}$, and
the $t$-sum is large. Later when we generalize the method we study
special curves that are quartic in $T$. Let
\begin{eqnarray}
y^2 \ = \  f(x,T) & \ = \  & x^3 T^2 + 2g(x)T - h(x) \nonumber\\
g(x) & \ = \  & x^3 + ax^2 + bx + c, \ c \neq 0 \nonumber\\ h(x) &
\ = \  & (A-1)x^3 + Bx^2 + Cx + D \nonumber\\ D_T(x) & \ = \  &
g(x)^2 + x^3 h(x).
\end{eqnarray}
Note that $D_T(x)$ is one-fourth of the discriminant of the
quadratic (in $T$) polynomial $f(x,T)$.
%The number of distinct, non-zero roots of $D_t(x)$ controls the rank.
 We write
$A-1$ as the leading coefficient of $h(x)$, and not $A$, to simplify future
computations by making the coefficient of $x^6$ in $D_T(x)$ equal
$A$.

Our elliptic curve $\mathcal{E}$ is not written in standard form,
as the coefficient of $x^3$ is $T^2$. This is harmless, and later
we rewrite the curve in Weierstrass form. As $y^2 = f(x,T)$, for
the fiber at $T=t$ we have
\begin{eqnarray}
a_t(p) \ = \  - \sum_{x (p)} \js{f(x,t)} \ = \  -\sum_{x (p)}
\js{x^3 t^2 + 2g(x)t - h(x)},
\end{eqnarray} where $\js{\ast}$ is the Legendre symbol.
We study $-pA_{\mathcal{E}}(p) = \zsum{x} \zsum{t} \js{f(x,t)}$.
When $x \equiv 0$ the $t$-sum vanishes if $c \not\equiv 0$, as it
is just $\zsum{t} \js{2ct - D}$. Assume now $x \not\equiv 0$. By
the lemma on Quadratic Legendre Sums (Lemma \ref{labquadlegsum})
\begin{equation}
\zsum{t} \js{x^3 t^2 + 2g(x)t - h(x)} \ = \  \Bigg\{
{(p-1)\js{x^3} \ \ \text{if} \ p \ | \ D_t(x) \atop \ \ -\js{x^3}
\ \ \text{otherwise} }
\end{equation}
Our goal is to find coefficients $a,b,c,A,B,C,D$ so that $D_t(x)$
has six distinct, non-zero roots. We want the roots $r_1, \dots,
r_6$ to be squares in $\Q$, as their contribution is $(p-1)\js{r_i^3}$. If
$r_i$ is not a square, $\js{r_i}$ will be $1$ for half the primes
and $-1$ for the other half, yielding no net contribution to the
rank. Thus, for $1 \leq i \leq 6$, let $r_i = \rho_i^2$.

Assume we can find such coefficients. Then
\begin{eqnarray}
-pA_{\mathcal{E}}(p) & \ = \   & \zsum{x} \zsum{t} \js{f(x,t)} \ =
\ \zsum{x} \zsum{t} \js{x^3 t^2 + 2g(x)t - h(x)} \nonumber\\ & \ =
\ & \sum_{x = 0} \zsum{t} \js{f(x,t)} + \sum_{x:D_t(x) \equiv 0}
\zsum{t} \js{f(x,t)} + \sum_{x:xD_t(x) \not\equiv 0} \zsum{t}
\js{f(x,t)} \nonumber\\ & \ = \  & 0 \ + \ 6(p-1) \ - \
\sum_{x:xD_t(x) \not\equiv 0} \js{x^3} \ = \  6p.
\end{eqnarray}
We must find $a, \dots, D$ such that $D_t(x)$ has six distinct,
non-zero roots $\rho_i^2$:
\begin{eqnarray}
D_t(x) & \ = \  & g(x)^2 + x^3 h(x) \nonumber\\ & \ = \  &  Ax^6 +
(B+2a) x^5 + (C + a^2 + 2b) x^4 + (D + 2ab + 2c) x^3 \nonumber\\ &
\ & + \ (2ac + b^2) x^2 + (2bc)x + c^2 \nonumber\\
%%Alvaro: I think the following line can be skept.
%%
%& \ = \  & A
%\Big( x^6 + \frac{B+2a}{A}x^5 + \frac{C+a^2+2b}{A}x^4 +
%\frac{D+2ab+2c}{A}x^3 \nonumber\\ & \ &  + \ \frac{2ac+b^2}{A}x^2
%+ \frac{2bc}{A}x + \frac{c^2}{A} \Big) \nonumber\\
& \ = \  &
A(x^6
+ R_5 x^5 + R_4 x^4 + R_3 x^3 + R_2 x^2 + R_1 x + R_0) \nonumber\\
& \ = \  &
A(x-\rho_1^2)(x-\rho_2^2)(x-\rho_3^2)(x-\rho_4^2)(x-\rho_5^2)(x-\rho_6^2).
\end{eqnarray}

\subsection{Determining Admissible Constants $a, \dots, D$}

Because of the freedom to choose $B,C,D$ there is no problem
matching coefficients for the $x^5, x^4, x^3$ terms. We must
simultaneously solve in integers
\begin{eqnarray}\label{eq:6abc}
2ac + b^2 & \ = \  & R_2 A \nonumber\\ 2bc & \ = \  & R_1 A \nonumber\\
c^2 & \ = \  & R_0 A.
\end{eqnarray} For simplicity, take $A = 64R_0^3$. Then
\begin{eqnarray}\label{eq:6abc}
\begin{array}{rrrrrrr}
c^2 & \ = \  & 64 R_0^4 & \longrightarrow  & c & \ = \  & 8 R_0^2  \\
2bc & \ = \  & 64R_0^3 R_1 & \longrightarrow  & b & \ = \  & 4 R_0 R_1  \\
2ac + b^2 & \ = \  & 64 R_0^3 R_2 & \longrightarrow  & a & \ = \ &
4R_0 R_2 - R_1^2.
\end{array}
\end{eqnarray}
For an explicit example, take $r_i = \rho_i^2 = i^2$. For these
choices of roots,
\begin{equation}
R_0 \ = \  518400,\ R_1 \ = \  -773136,\ R_2 \ = \  296296.
\end{equation}
Solving for $a$ through $D$ yields
\begin{equation}\label{eq:6all}
\begin{array}{ccrrrrrr}
A & \ = \  & 64 R_0^3 & \ = \  & 8916100448256000000 & \\
c & \ = \  & 8 R_0^2 & \ = \  & 2149908480000 & \\ b & \ = \ & 4
R_0 R_1 & \ = \  & -1603174809600 & \\ a & \ = \ & 4 R_0 R_2 -
R_1^2 & \ = \  & 16660111104 & \\ B & \ = \  & R_5 A - 2a & \ = \
& \ -811365140824616222208 & \\ C & \ = \  & R_4 A - a^2 - 2b & \
= \  & \ 26497490347321493520384 & \\ D & \ = \  & R_3 A - 2ab -
2c & \ = \  & -343107594345448813363200 &
\end{array}
\end{equation}

We convert $y^2 = f(x,t)$ to $y^2 = F(x,t),$ which is in
Weierstrass normal form. We send $y \to \frac{y}{t^2 + 2t - A +
1},$ $x \to \frac{x}{t^2 + 2t - A + 1},$ and then multiply both
sides by $(t^2 + 2t - A + 1)^2$. For future reference, we note
that
\begin{eqnarray}
t^2 + 2t - A + 1 & \ = \  & (t + 1 - \sqrt{A})(t + 1 + \sqrt{A})
\nonumber\\ & \ = \  & (t - t_1)(t - t_2) \nonumber\\ & \ = \  &
(t - 2985983999)(t + 2985984001).
\end{eqnarray}
We have
\begin{eqnarray}
& f(x,t) & \ = \  t^2 x^3 + (2x^3 + 2ax^2 + 2bx + 2c)t - (A-1)x^3
- Bx^2 - Cx - D \nonumber\\ & \ & \ = \  (t^2 + 2t - A + 1)x^3 +
(2at - B)x^2 + (2bt - C)x + (2ct - D) \nonumber\\ & F(x,t) & \ = \
x^3 + (2at-B)x^2 + (2bt-C)(t^2 + 2t - A + 1)x \nonumber\\ & \ & \
\ \ \ \ \ + (2ct - D)(t^2 + 2t - A + 1)^2.
\end{eqnarray}
We now study the $-pA_{\mathcal{E}}(p)$ arising from $y^2 =
F(x,T)$. It is enough to show this is $6p + O(1)$ for all $p$
greater than some $p_0$. Note that $t_1, t_2$ are the unique roots
of $t^2+2t-A+1 \equiv 0\bmod p$. We find \bea -pA_{\mathcal{E}}(p)
\ = \ \zsum{t} \zsum{x} \js{F(x,t)} & \ = \  & \sum_{t \neq t_1,
t_2} \zsum{x} \js{F(x,t)} + \sum_{t = t_1,t_2} \zsum{x}
\js{F(x,t)}. \eea For $t \neq t_1,t_2$, send $x \longrightarrow
(t^2 + 2t - A + 1)x$. As $(t^2 + 2t - A + 1) \not\equiv 0$,
$\js{(t^2 + 2t - A + 1)^2} = 1$. Simple algebra yields \bea
-pA_{\mathcal{E}}(p) & \ = \ & 6p +
O(1) + \sum_{t=t_1,t_2} \zsum{x} \js{f_{t}(x)} + O(1) \nonumber\\
& \ = \ & 6p + O(1) + \sum_{t=t_1,t_2} \zsum{x} \js{(2at-B)x^2 +
(2bt-C)x + (2ct-D)}.
\end{eqnarray}
The last sum above is negligible (i.e., is $O(1)$) if
\begin{equation}\label{eq:tcheck}
D(t) \ = \  (2bt - C)^2 - 4(2at - B)(2ct - D) \not\equiv 0 (p).
\end{equation}
Calculating yields
\begin{eqnarray}
D(t_1) & \ = \  & 4291243480243836561123092143580209905401856
\nonumber\\ &\ = \ & 2^{32}\cdot 3^{25} \cdot7^5 \cdot 11^2\cdot
13 \cdot 19 \cdot 29 \cdot
31 \cdot 47 \cdot 67 \cdot 83 \cdot 97 \cdot 103 \nonumber\\
D(t_2) & \ = \  & 4291243816662452751895093255391719515488256
\nonumber\\ &\ = \ & 2^{33}\cdot 3^{12}\cdot 7 \cdot 11 \cdot 13
\cdot 41 \cdot 173 \cdot 17389 \cdot 805873 \cdot 9447850813. \ \
\
\end{eqnarray}
Hence, except for finitely many primes (coming from factors of
$D(t_i)$, $a, \dots, D$, $t_1$ and $t_2$), $-A_\mathcal{E}(p) = 6p
+ O(1)$ as desired. We have shown the following result:
\begin{thm}
There exist integers $a,b,c,A,B,$ $C,D$ so that the curve
$\mathcal{E}:$ $y^2 = x^3 T^2 + 2g(x)T - h(x)$ over $\Q(T)$, with
$g(x) = x^3 + ax^2 + bx + c$ and $h(x) = (A-1)x^3 + Bx^2 + Cx +
D$, has rank $6$ over $\Qoft$. In particular, with the choices of
$a$ through $D$ above, $\mathcal{E}$ is a rational elliptic
surface and has Weierstrass form
$$y^2 = x^3 + (2aT-B)x^2 + (2bT-C)(T^2 + 2T - A + 1)x
 + (2cT - D)(T^2 + 2T - A + 1)^2$$
%%Alvaro: DO WE NEED ALL THESE QUANTITIES BELOW IN THE PAPER?
%c_4(T) &\ = \ & 2^{19}3^7 7^1 13^1 (1475 T^3 +
%12359745382011 T^2 \ - \ 4860110603997053240403 T
%\nonumber\\ & & \ - \ 7735999878503076170786750620939) \nonumber\\
%c_6(T) & \ = \  & -2^{25}3^{11}(625 T^5 + \cdots ) \nonumber\\
%j(T) & \ = \  & \frac{50141357421875 T^{9} +
%\cdots}{-1171875T^{10} + \cdots } \nonumber\\
%\Delta(T) & \ = \  & -2^{44}3^{18} 5^6(75 T^{10} + \cdots) \ = \
%-2^{44}3^{18} 5^6 (T \pm \sqrt{A} + 1)^2(75T^6 + \cdots ).
%\end{eqnarray}
\end{thm}

\begin{proof} We show $\mathcal{E}$ is a rational
elliptic surface by translating $x\mapsto x-(2aT-B)/3$, which
yields $y^2=x^3+A(T)x+B(T)$ with $\deg(A)=3, \deg(B)=5$. Therefore
the Rosen-Silverman theorem is applicable, and because we can
compute $A_\mathcal{E}(p)$, we know the rank is exactly 6 (and we
never need to calculate height matrices).
\end{proof}

\begin{rek} We can construct infinitely many $\mathcal{E}$ over
$\Q(T)$ with rank 6 using \eqref{eq:6all}, as for generic choices
of roots $\rho_1^2, \dots, \rho_6^2$, \eqref{eq:tcheck} holds.
\end{rek}

For concreteness, we explicitly list a curve of rank at least $6$.
Doing a better job of choosing coefficients $a$ through $D$ (but
still being crude) yields

\begin{thm}
$y^2 = x^3 + Ax + B$ has rank at least 6, where
\begin{eqnarray}
A &  \ = \  & \ \ \ \ \ \ \ \ \ \ \ \ \ \ 1123187040185717205972
\nonumber\\ B & \ = \  & 50786893859117937639786031372848
\nonumber\
\end{eqnarray}
\end{thm}

Six points on the curve are: \be \begin{array}{rrrr}
       ( 67585071288,\ & 20866449849961716) & ( 60673071396,\ & 18500949214922664)  \\
       ( 49153071576,\ & 14991664661755236) & ( 33025071828,\ & 11131001682078096)   \\
       ( 12289072152,\  & 8151425152633980) & (-13054927452,\ & 5822267813027064)
\end{array} \ee

As the determinant of the height matrix is approximately
$880,000$, the points are independent and therefore generate the
group. A trivial modification of this procedure yields rational
elliptic surfaces of any rank $r \le 6$. For more constructions
along these lines, see \cite{Mil}.

%!!!!!!!!!!!!!!!!!!!!!!!!!!!!!!!!!!!!!!!new section!!!!!!!!!!!!!!!
%!!!!!!!!!!!!!!!!!!!!!!!!!!!!!!!!!!!!!!!new section!!!!!!!!!!!!!!!
%!!!!!!!!!!!!!!!!!!!!!!!!!!!!!!!!!!!!!!!new section!!!!!!!!!!!!!!!

\section{More Attempts for Curves with rank $6$, $7$ and $8$ over $\Q(T)$}
\setcounter{equation}{0}

\subsection{Curves of Rank $6$}

We sketch another construction for a curve of rank $6$ over
$\Q(T)$ by modifying our previous arguments. We define a curve
$\mathcal{E}$ over $\Q(T)$ by
\begin{eqnarray}
y^2 \ = \  f(x,T) & \ = \  & x^4 T^2 + 2g(x)T - h(x) \nonumber\\
g(x) & \ = \  & x^4 + ax^3 + bx^2 + cx + d, \ d \neq 0 \nonumber\\
h(x) & \ = \  & -x^4 + Ax^3 + Bx^2 + Cx + D \nonumber\\ D_T(x) & \
= \ & g(x)^2 + x^4 h(x).
\end{eqnarray}

We must find choices of the free coefficients such that $D_t(x) =
\prod_{i=1}^7 (\ga^2 x-\rho_i)$, with each root non-zero. For
$x=0$, we have $\sum_t \js{2dt - D} = 0$. By Lemma
\ref{labquadlegsum}, for $x$ a root of $D_t$ we have a
contribution of $(p-1)\js{x^4} = (p-1)\js{\rho_i^4 \alpha^{-8} } =
p-1$; for all other $x$ a contribution of $-\js{x^4\alpha^{-8}} =
-1$. Hence summing over $x$ and $t$ yields $7(p-1) + \sum_{x \neq
\rho_i, 0} -1 = 6p$. Similar reasoning as before shows we can find
integer solutions (we included the factor of $\ga^2$ to facilitate
finding such solutions). We  chose the coefficient  of the $x^4$
term to be $t^2 + 2t + 1 = (t+1)^2$, as this implies each curve
$E_t$ is isomorphic over $\Q$ to an elliptic curve $E_t'$ (see
Appendix \ref{sec:convertqc}). As $\mathcal{E}$ is almost
certainly not rational, the rank is exactly 6 if Tate's conjecture
is true for the surface. If we only desire a lower bound for the
rank, we can
 list the 6 points and calculate the determinant of the height
matrix and see if they are independent.

\subsection{Probable Rank $7$, $8$ Curves}

We modify the previous construction to
\begin{eqnarray}
y^2 & \ = \  & x^3 T^2 + 2g(x) T - h(x) \nonumber\\ g(x) & \ = \ &
x^4 + a x^3 + b x^2 + cx + d, \ d \neq 0 \nonumber\\ h(x) & \ = \
& Ax^4 + Bx^3 + Cx^2 + Dx + E
\end{eqnarray}
to obtain what should be higher rank curves over $\Q(T)$. Choosing
appropriate quartics for $g(x), h(x)$ such that $D_T(x) = g^2(x) +
x^3 h(x)$ has eight distinct non-zero perfect square roots should
yield a contribution of $8p$. As the coefficient of $T^2$ is
$x^3$, we do not lose $p$ from summing over non-roots of $D_t(x)$.
By specializing to certain $t = a_2 s^2 + a_1 s + a_0$ for some
constants, we can arrange it so $y^2 = k^2(s) x^4 + \cdots$, and
by the previous arguments obtain a cubic. Unfortunately, we can no
longer explicitly evaluate $pA_{\mathcal{E}}(p)$ (because of the
replacement $t \to a_2 s^2 + a_1 s + a_0$). As the method yields
eight points for all $s$, we need only specialize and compute the
height matrix. As we construct a rank 8 curve over $\Q(T)$ in
\S\ref{sec:usingquartics} (when we generalize our construction),
we do not provide the details here. Note, however, that sometimes
there are obstructions and the rank is lower than one would expect
(see \S\ref{sec:lindep}).

\section{Using Cubics and Quartics in $T$}\label{sec:usingquartics}

Previously we used $y^2 = f(x,T)$, with $f$ quadratic in $T$. The
reason is that, for special $x$, we obtain $y_i^2 = s_i(x_i)^2 (T
- t_i)^2$. For such $x$, the $t$-sum is large (of size $p$); we
then show for other $x$ that the $t$-sum is small.

\subsection{Idea of Construction}

The natural generalization of our Discriminant Method is to
consider $y^2 = f(x,T)$, with $f$ of higher order in $T$. We first consider
polynomials cubic in $T$. For a fixed $x_i$,
we have the $t$-sum $\sum_{t (p)} \js{f(x_i,t)}$, and there are
several possibilities:

\ben \item $f(x_i,T) = a (T - t_1)^3$. In this case, the $t$-sum
will vanish, as $\js{ (t-t_1)^3} = \js{ t-t_1}$. \item $f(x_i,T) =
a (T - t_1)^2 (T - t_2)$. The $t$-sum will be $O(1)$, as for $t
\neq t_1$ we have $\js{ (t-t_1)^2 (t-t_2)} = \js{ t-t_2}$. \item
$f(x_i,T) = a(T-t_1) (T-t_2) (T - t_3)$. This will in general be
of size $\sqrt{p}$. \item $f(x_i,T) = a (T - t_1) (T^2 + bT + c)$,
with the quadratic irreducible over $\Z/p\Z$. This happens when
$b^2 - 4c$ is not a square mod $p$. This will in general be of
size $\sqrt{p}$. \item $f(x_i,T) = a T^3 + b T^2 + c T + d$, with
the cubic irreducible over $\Z/p\Z$. Again, this will in general
be of size $\sqrt{p}$. \een

Thus, our method does not generalize to $f(x,T)$ cubic in $T$. The
problem is we cannot reduce to $\js{ (t - t_1)^{2n_1} \cdots (t -
t_i)^{2n_i}}$. We therefore investigate $f(x,T)$ quartic in $T$.
Consider, for simplicity, a curve $\mathcal{E}$ over $\Q(T)$ of
the form: \be y^2 \ = \ f(x,T) \ = \ A(x) T^4 + B(x) T^2 + C(x),
\ee $A(x), B(x), C(x) \in \Z[x]$ of degree at most $4$. The
polynomial $A T^4 + B T^2 + C$ has discriminant $16 AC (4 A C -
B^2 )^2$. There are several possibilities for special choices of
$x$ giving rise to large $t$-sums (sums of size $p$): \ben \item
$A(x_i), B(x_i) \equiv 0 \bmod p$, $C(x_i)$ a non-zero square mod
$p$. Then the $t$-summand is of the form $c^2$, contributing $p$.
\item $A(x_i), C(x_i) \equiv 0  \bmod p$, $B(x_i)$ a non-zero
square mod $p$. Then the $t$-summand is of the form $(bt)^2$,
contributing $p-1$. \item $B(x_i), C(x_i) \equiv 0 \bmod p$,
$A(x_i)$ a non-zero square mod $p$. Then the $t$-summand is of the
form $(at^2)^2$, contributing $p-1$. \item $A(x_i)$ is a non-zero
square mod $p$ and $B(x_i)^2 - 4 A(x_i) C(x_i) \equiv 0 \bmod p$.
Then the $t$-summand is of the form $a^2 (t^2 - t_1)^2$,
contributing $p-1$. \een

In the above construction, we are no longer able to calculate
$A_{\mathcal{E}}(p)$ exactly. Instead, we construct curves where
we believe $A_{\mathcal{E}}(p)$ is large. This is accomplished by
forcing points to be on $\mathcal{E}$ which satisfy any of (1)
through (4) above. As we are unable to evaluate the
$A_\mathcal{E}(p)$ sums, we specialize and calculate height
matrices to show the points are independent. Unfortunately, some
of our constructions yielded 9 and 10 points on $\mathcal{E}$, but
some of these points were linearly dependent on the others, or
torsion points (see \S\ref{sec:lindep}).

This method, with a quartic in $T$, can force a maximum number of
12 points on $\mathcal{E}$. It is possible to have 8 points from
the vanishing of the discriminant (in $t$), and an additional 6
points from the simultaneous vanishing of pairs of $A(x), B(x),
C(x)$; however, any common root of $A$ or $C$ with $B$ is also a
root of $B^2-4AC$, so there are at most 4 new roots arising from
simultaneous vanishing, for a total of 12 possible points.

\subsection{Rank (at least) 7 Curve}

For appropriate choices of the parameters, the curve $\mathcal{E}:
y^2 \ = \ A(x)T^4 + 4 B(x)T^2 + 4 C(x)$ over $\Q(T)$ with
\begin{eqnarray}
A(x)&=&a_1a_2a_3a_4(x-a_1)(x-a_2)(x-a_3)(x-a_4)\nonumber\\
C(x)&=&a_1a_2c_1c_2(x-a_1)(x-a_2)(x-c_1)(x-c_2)\nonumber\\
B(x)&=&a_1^2a_2^2(x-c_1)(x-c_2)(x-a_3)(x-a_4)
\end{eqnarray}
has rank at least 7. We get 6 points from the common vanishing of
$A,B,C$ in pairs and an additional point from a factor of
$B^2-AC$. Choosing $a_1=-25, a_2=-5, a_3=-10, a_4=-1, c_1=-9,
c_2=15$ we find that the points
\begin{eqnarray}
(-25,120000T), (-5,10000T), (-10,11250), (-1,28800),\nonumber\\
(-9,800T^2), (15, 20000T^2), (65/7, (540000t^2- 2880000)/49)
\end{eqnarray}
all lie on $\mathcal{E}$. Upon transforming to a cubic (see
Appendix \ref{sec:convertqc}), specializing to $T=20$, and
considering the minimal model, we found that these points are
linearly independent (PARI calculates the determinant of the
height matrix  is approximately $37472$). Note this is not a
rational surface, as the coefficient of $x$ in Weierstrass form is
of degree 8.

\subsection{Rank (at least) 8 Curve}

For appropriate choices of the parameters, the curve $\mathcal{E}:
y^2 \ = \ A(x)T^4 + B(x)T^2+C(x)$ over $\Q(T)$ with  \be \nonumber
A(x)=x^4,\ B(x)=2x(b_3x^3+b_2x^2+b_1x+b_0)+b^2,\ C(x)=
x(b_3^2x^3+c_2x^2+c_1x+c_0) \ee has rank at least 8. As the
coefficient of $x^4$ is $T^4+2b_3T^2+b_3^2$, a perfect square,
$\mathcal{E}$ can easily be transformed into Weierstrass form (see
Appendix \ref{sec:convertqc}). The common vanishing of $A$ and $C$
at $x=0$ produces a point $S_0=(0,bT)$ on
$\mathcal{E}/\mathbb{Q}(T)$. Also notice that as before, if
$B^2-4AC$ vanishes at $x=x_i$ then we can rewrite: \be
A(x_i)T^4+B(x_i)T^2+C(x_i)\ = \
A(x_i)\left(T^2+\frac{B(x_i)}{2A(x_i)}\right)^2\ = \
x_i^4\left(T^2+\frac{B(x_i)}{2x_i^4}\right)^2\ee Thus we obtain a
point $P_{x_i}=(x_i,x_i^2(T^2 +B(x_i)/2x_i^4))$ on $\mathcal{E}$.
We chose constants $b_i,b$ an $c_i$ so that  \be
B^2-4AC=(x-1)(x+1)(x-4)(x+4)(x-9)(x+9)(x-16)\ee and obtain  a
curve $\mathcal{E}$ over $\Q(T)$ with coefficients: \bea
\label{rank8} A & \ = \ & x^4, \quad B(x) =
-\frac{5852770213}{382205952}x^4 + \frac{89071}{36864}x^3 -
\frac{89233}{1152}x^2 - \frac{9}{2}x + 144, \nonumber\\  C(x) & =
& \frac{34254919166180065369}{584325558976905216}x^4 -
\frac{528356915749387}{28179280429056}x^3 \nonumber\\\ & & +
\frac{527067904642903}{880602513408}x^2 -
\frac{5881576729}{169869312}x.  \eea

As discussed above, the curve $\mathcal{E}$ given by (\ref{rank8})
has $8$ rational points over $\Q(T)$, namely $S_0$ and $P_{x_i}$
for $x_i=\pm 1,\ \pm 4,\ \pm 9, 16$. As $\mathcal{E}$ is not  a
rational surface, and as we cannot evaluate $A_{\mathcal{E}}(p)$
exactly, we need to make sure the points are linearly independent.
Specializing to $T = 1$ yields the elliptic curve with minimal
model
\be \begin{array}{rcr} E_1\colon y^2 &=& x^3-x^2-\alpha x+\beta \\
\alpha &=& 357917711928106838175050781865\\
\beta &=& 8790806811671574287759992288018136706011725.\end{array}
\ee The eight points of $E_T$ at $T=1$ are linearly independent on
$E_1/\mathbb{Q}$ (PARI calculates the determinant of the height
matrix to be about $124079248627.08$), proving $\mathcal{E}$ does
have rank at least 8 over $\Q(T)$.

\section{Linear Dependencies Among Points}\label{sec:lindep}

Not all choices of $A(x), B(x), C(x)$ which yield $r$ points on
the curve $\mathcal{E}: y^2 = A(x)T^4 + 4 B(x)T^2 + 4 C(x)$
actually give a curve of rank at least $r$ over $\Q(T)$. We found
many examples giving 9 and 10 points by choosing $A(x) = C(x)$ so
that $B^2-AC$ factors nicely, and then searching through
prospective roots of this quantity as well as roots of $A(x) =
C(x)$. One such curve giving 10 points arises from
\begin{eqnarray}
A(x)&\ = \ &C(x)\ = \ (x-1)^2(2x-1)^2\nonumber\\
B(x)&=&12316x^4+2346x^3-239x^2-24x+1,
\end{eqnarray}
and has the following points on it
\begin{eqnarray}
(0, T^2+2), \left( \frac{-1}{19}, \frac{420}{361}(T^2+2)\right),
\left( \frac{-1}{4}, \frac{15}{8}(T^2+2)\right),\nonumber\\
\left( \frac{1}{9}, \frac{56}{81}(T^2+2)\right), \left(
\frac{-1}{7},  \frac{72}{49}(T^2-2)\right), \left( \frac{-1}{5},
 \frac{42}{25}(T^2-2)\right),\nonumber\\
\left( \frac{1}{11}, \frac{90}{121}(T^2-2)\right), \left(
\frac{1}{16},  \frac{105}{128}(T^2-2)\right), (1, 240T), \left(
\frac{1}{2}, 63T\right).
\end{eqnarray}

It can be shown, however, that upon translating to a cubic only
the (translated versions of the) second, third, fifth, sixth, and
ninth of these points are independent over $\mathbb{Q}(T)$. While
the contribution from these points makes $A_{\mathcal{E}}(p)$ want
to be large, this is not reflected by a large rank.

\section{Using Higher Degree Polynomials}

Let $f(x,T)$ be a polynomial of degree 3 or 4 in $x$ and arbitrary
degree in $T$ and let $\mathcal{E}$ be the elliptic curve over
$\Q(T)$ given by $y^2=f(x,T)$ (with the coefficient of $x^4$ a
perfect square or zero). The remarks at the beginning of Section
\ref{sec:usingquartics} about cubics suggest that we should look
for polynomials $f(x,T)$ with even degree in $T$, say
$\deg_T(f)=2n$.

The nice feature of quadratics and biquadratics that we used in
the previous constructions was the fact that a zero of the
discriminant indicates that the polynomial $f(x,T)$ factors as a
perfect square. However, when $f$ is of arbitrary degree $2n$ in
$T$ this is no longer true: a zero of the discriminant $D_T(x)$
indicates just a multiple root. However, in the most general case,
there exist $n$ quantities $D_{i,T}(x)$ such that their common
vanishing at $x=x_0$ implies that $f(x,T)$ factors as a perfect
square. As an example we look at a quartic of the form
$f(x,T)=A^2T^4+BT^3+CT^2+DT+E^2$, where $\deg_x(A,E)\leq 2$ and
$\deg_x(B,C,D)\leq 4$. This can be rewritten as:

\begin{center}
$ A^2T^4 + 2AT^2( \frac{Bt}{2A} + \frac{C}{2A} - \frac{B^2}{8A^3}
) + ( \frac{BT}{2A} + \frac{C}{2A} - \frac{B^2}{8A^3} )^2 + (D -
\frac{B}{A}( \frac{C}{2A} - \frac{B^2}{8A^3} ))T- (\frac{C}{2A}
-\frac{B^2}{8A^3})^2 + E^2.$
\end{center}
The last two terms are the ones which are keeping the polynomial
from being a perfect square. Thus, if \be D - \frac{B}{A}\left(
\frac{C}{2A} - \frac{B^2}{8A^3} \right)\ =\ 0,\quad E^2 -
\left(\frac{C}{2A} - \frac{B^2}{8A^3}\right)^2\ = \ 0 \ee then the
polynomial $f$ will be a square. This is equivalent to
\begin{eqnarray}
D_{1,T} & = & 8A^4 D-4A^2 BC+B^3 = 0\nonumber\\
D_{2,T} & = & 64A^6 E^2 - 16A^4 C^2 - B^4 + 8A^2 C B^2 =0.
\end{eqnarray}
Note that if B=D=0, the conditions that these polynomials impose
reduce to the usual discriminant. Also, $\deg_x(D_{1,T})\leq 12,\
\deg_x(D_{2,T})\leq 16$, so we could get up to $12$ points of
common vanishing of the $D_i$. The authors have tried to find
suitable constants without success, due to the complexity of the
Diophantine equations.

\section*{Acknowledgements}

We thank Vitaly Bergelson, Michael Hunt, James Mailhot, David
Rohrlich, Mustafa Sahin, and Warren Sinnott for many enlightening
conversations. The second author also wishes to thank Boston
University for its hospitality, where much of the final write-up
was prepared.

%!!!!!!!!!!!!!!!!!!!!!!!!!!!!!!!!!!!!!!!new section!!!!!!!!!!!!!!!!!!!!!!!!!!!!!!!!
%!!!!!!!!!!!!!!!!!!!!!!!!!!!!!!!!!!!!!!!new section!!!!!!!!!!!!!!!!!!!!!!!!!!!!!!!!
%!!!!!!!!!!!!!!!!!!!!!!!!!!!!!!!!!!!!!!!new section!!!!!!!!!!!!!!!!!!!!!!!!!!!!!!!!

\appendix

\section{Sums of Legendre Symbols}

For completeness, we provide proofs of the quadratic Legendre sums
that are used in our constructions.

\subsection{Factorizable Quadratics in Sums of Legendre Symbols}

\begin{lem}\label{lablegsum} For $p > 2$
\begin{equation}
\twocase{S(n) \ = \  \zsum{x} \js{n_1 + x} \js{n_2 + x} \ = \
}{p-1}{if $p | (n_1 - n_2)$}{-1}{otherwise.}
\end{equation}
\end{lem}
\begin{proof} Translating $x$ by $-n_2$, we need only prove the lemma
when $n_2 =   0$. Assume $(n,p) = 1$ as otherwise the result is
trivial. For $(a,p) = 1$ we have:
\begin{eqnarray}
S(n) & \ = \  & \sum_{x=0}^{p-1} \js{n + x} \js{x} \nonumber\\
     & \ = \  & \zsum{x} \js{n + a^{-1} x} \js{a^{-1} x} \nonumber\\
     & \ = \  & \zsum{x} \js{an + x} \js{x} \ = \  S(an)
\end{eqnarray}
Hence
\begin{eqnarray}
S(n) & \ = \  & \frac{1}{p-1} \osum{a} \zsum{x} \js{an+x} \js{x} \nonumber\\
     & \ = \  & \frac{1}{p-1} \zsum{a} \zsum{x} \js{an+x} \js{x} -
     \frac{1}{p-1} \zsum{x} \js{x}^2 \nonumber\\
     & \ = \  & \frac{1}{p-1} \zsum{x} \js{x} \zsum{a} \js{an+x} - 1
     \nonumber\\
     & \ = \  & 0 - 1 \ = \  -1.\end{eqnarray}\end{proof}

We need $p > 2$ as we used $\sum_{a=0}^{p-1} \js{an+x} = 0$ for
$(n,p) = 1$. This is true for all odd primes (as there are
$\frac{p-1}{2}$ quadratic residues, $\frac{p-1}{2}$ non-residues,
and $0$); for $p=2$, there is one quadratic residue, no
non-residues, and $0$.

\subsection{General Quadratics in Sums of Legendre Symbols}

\begin{lem}[Quadratic Legendre Sums]\label{labquadlegsum} Assume
$a$ and $b$ are not both zero mod $p$ and $p > 2$. Then
\begin{equation}
\twocase{\zsum{t} \js{at^2 + bt + c} \ = \ }{(p-1)\js{a}}{if $p |
 (b^2 - 4ac)$}{-\js{a}}{otherwise.}
\end{equation}
\end{lem}

\begin{proof} Assume $a \not\equiv 0 (p)$ as otherwise the proof is
trivial. Let $\delta \ = \  4^{-1}(b^2 - 4ac)$. Then
\begin{eqnarray}
\zsum{t} \js{at^2 + bt + c} & \ = \  & \zsum{t} \js{a^{-1}}
\js{a^2 t^2 + bat + ac} \nonumber\\ & \ = \  & \zsum{t} \js{a}
\js{t^2 + bt + ac} \nonumber\\  & \ = \  & \zsum{t} \js{a} \js{ (t
+ 2^{-1}b)^2 - 4^{-1}(b^2-4ac)} \nonumber\\ & = & \js{a} \zsum{t}
\js{t^2 -\delta}
\end{eqnarray}
If $\delta \equiv 0 (p)$ we get $p-1$. If $\delta \equiv \eta^2,
\eta \neq 0$, then by Lemma \ref{lablegsum}
\begin{equation}
\zsum{t} \js{t^2 - \delta} \ = \  \zsum{t} \js{t-\eta} \js{t+\eta}
\ = \  -1.
\end{equation}
We note that $\zsum{t} \js{t^2 - \delta}$ is the same for all
non-square $\delta$'s (let $g$ be a generator of the
multiplicative group, $\delta = g^{2k+1}$, change variables by $t
\to g^k t$). Denote this sum by $S$, the set of non-zero squares
mod $p$ by $\mathcal{R}$, and the non-squares mod $p$ by
$\mathcal{N}$. Since $\zsum{\delta} \js{t^2 - \delta} = 0$ we have
\begin{eqnarray}
\zsum{\delta} \zsum{t} \js{t^2 - \delta} & \ = \  & \zsum{t}
\js{t^2} + \sum_{\delta \in \mathcal{R}} \zsum{t} \js{t^2 -
\delta} + \sum_{\delta \in \mathcal{N}} \zsum{t} \js{t^2 - \delta}
\nonumber\\ & \ = \  & (p-1) + \frac{p-1}{2}(-1) + \frac{p-1}{2}S
\ = \  0
\end{eqnarray}
Hence $S = -1$, proving the lemma.\end{proof}

%!!!!!!!!!!!!!!!!!!!!!!!!!!!!!!!!!!!!!!!new section!!!!!!!!!!!!!!!!!!!!!!!!!!!!!!!!
%!!!!!!!!!!!!!!!!!!!!!!!!!!!!!!!!!!!!!!!new section!!!!!!!!!!!!!!!!!!!!!!!!!!!!!!!!
%!!!!!!!!!!!!!!!!!!!!!!!!!!!!!!!!!!!!!!!new section!!!!!!!!!!!!!!!!!!!!!!!!!!!!!!!!

\section{Converting from Quartics to Cubics}\label{sec:convertqc}

We record two useful transformations from quartics to cubics. In
all theorems below, all quantities are rational.

\begin{thm} If the quartic curve $y^2 = x^4 - 6cx^2 + 4dx + e$
has a rational point, then it is equivalent to the cubic curve
$Y^2 = 4X^3 - g_2X - g_3$, where \be g_2\ =\ e+3c^2, \ \ \ \ g_3\
= \ -ce - d^2 + c^3, \ee and \be 2x \ = \ (Y-d)/(X-c), \ \ \ \ y \
= \ -x^2 + 2X + c.\ee \end{thm}

See \cite{Mor}, page 77. Note that if the leading term of the
quartic is $a^2x^4$, one can send $y\to y/a$ and $x\to x/a$.

\begin{thm} The quartic $v^2 = au^4 + bu^3 + cu^2 + du + q^2$ is
equivalent to the cubic $y^2 + a_1xy + a_3y = x^3 +
a_2x^2+a_4x+a_6$, where \be a_1 = d/q, \ \ a_2 = c - (d^2/4q^2), \
\ a_3 = 2qb, \ \ a_4 = -4q^2a, \ \ a_6 = a_2a_4 \ee and \be x \ =
\ \frac{2q(v+q) + du}{u^2}, \ \ \ \ y\ = \
\frac{4q^2(v+q)+2q(du+cu^2)-(d^2u^2/2q)}{u^3}. \ee The point
$(u,v)=(0,q)$ corresponds to $(x,y)=\infty$ and $(u,v)=(0,-q)$
corresponds to $(x,y)=(-a_2,a_1a_2-a_3)$.
\end{thm}

See \cite{Wa}, page 37.

%!!!!!!!!!!!!!!!!!!!!!!!!!!!!!!!!!!!!!!!new section!!!!!!!!!!!!!!!!!!!!!!!!!!!!!!!!
%!!!!!!!!!!!!!!!!!!!!!!!!!!!!!!!!!!!!!!!new section!!!!!!!!!!!!!!!!!!!!!!!!!!!!!!!!
%!!!!!!!!!!!!!!!!!!!!!!!!!!!!!!!!!!!!!!!new section!!!!!!!!!!!!!!!!!!!!!!!!!!!!!!!!


\begin{thebibliography}{PTW02} % '2nd argument contains the widest acronym'

\bibitem[Fe1]{Fe1}
\newblock S. Fermigier, \emph{\'Etude exp\'erimentale du rang de
familles de courbes elliptiques sur $\Q$}, Exper. Math.
\textbf{5}, $1996$, $119-130$.

\bibitem[Fe2]{Fe2}
\newblock S. Fermigier, \emph{Une courbe elliptique d\'{e}finie
sur $\Q$ de rang $\ge 22$}, Acta Arith. \textbf{82} 1997, no. 4,
359–363.

\bibitem[Mes1]{Mes1}
\newblock J. Mestre, \emph{Courbes elliptiques de rang $\ge$ 11 sur $\Q(T)$},
C. R. Acad. Sci. Paris, ser. $1$, \textbf{313}, $1991$, $139-142$.

\bibitem[Mes2]{Mes2}
\newblock J. Mestre, \emph{Courbes elliptiques de rang $\ge$ 12 sur $\Q(T)$},
C. R. Acad. Sci. Paris, ser. $1$, \textbf{313}, $1991$, $171-174$.

\bibitem[Mil]{Mil}
\newblock S. J. Miller, \emph{$1$- and $2$-Level Densities for Families of Elliptic
Curves: Evidence for the Underlying Group Symmetries}, Ph.D.
Thesis, Princeton University, $2002$,
http://www.math.princeton.edu/$\sim$sjmiller/thesis/thesis.ps.

\bibitem[Mir]{Mir}
\newblock R. Miranda, \emph{The Basic Theory of Elliptic
Surfaces}, Dottorato di Ricerca in Matematica, Dipartmento di
Matematica dell'Universit'a di Pisa, ETS Editrice, $1989$.

\bibitem[Mor]{Mor}
\newblock Mordell, \emph{Diophantine Equations}, Academic Press, New York, $1969$,

\bibitem[Na1]{Na1}
\newblock K. Nagao, \emph{Construction of high-rank elliptic curves},
 Kobe J. Math. \textbf{11}, $1994$, $211-219$.

\bibitem[Na2]{Na2}
\newblock K. Nagao, \emph{$\Qoft$-rank of elliptic curves and certain
limit coming from the local points}, Manuscr. Math. \textbf{92},
$1997$, $13-32$.

%\bibitem[Nag]{Nag}
%\newblock T. Nagell, \emph{Introduction to Number Theory}, Chelsea
%Publishing Company, New York, $1981$.

\bibitem[RS]{RS}
\newblock M. Rosen and J. Silverman, \emph{On the rank of an elliptic
surface}, Invent. Math. \textbf{133}, $1998$, $43-67$.

\bibitem[RuS]{RuS}
\newblock K. Rubin and A. Silverberg, \emph{Ranks of elliptic
curves}, Bulletin of the American Mathematical Society
\textbf{39}, number $4$, $2002$, $455-474$.

\bibitem[Sh]{Sh}
\newblock T. Shioda, \emph{Construction of elliptic curves with high-rank
via the invariants of the Weyl groups}, J. Math. Soc. Japan
\textbf{43}, $1991$, $673-719$.

\bibitem[Si1]{Si1}
\newblock J. Silverman, \emph{The Arithmetic of
Elliptic Curves}, Graduate Texts in Mathematics \textbf{106},
Springer-Verlag, Berlin - New York, $1986$.

\bibitem[Si2]{Si2}
\newblock J. Silverman, \emph{Advanced Topics in the Arithmetic of
Elliptic Curves}, Graduate Texts in Mathematics \textbf{151},
Springer-Verlag, Berlin - New York, $1994$.

%\bibitem[Si3]{Si3}
%\newblock J. Silverman, \emph{The average rank of an algebraic
%family of elliptic curves}, J. reine angew. Math. \textbf{504},
%$1998$, $227-236$.

\bibitem[Wa]{Wa}
\newblock L. Washington, \emph{Elliptic Curves: Number Theory and Cryptography
(Discrete Mathematics and Its Applications)}, Chapman \& Hall,
2003.



\end{thebibliography}
\end{document}